\documentclass[12pt]{amsart}
\usepackage{amscd,amssymb,graphics}
\usepackage[matrix,arrow,curve]{xy}
\newcommand\mathscr{\mathcal}
\newtheorem{theorem}{Theorem}[section]
\newtheorem{corollary}[theorem]{Corollary}

\newtheorem{proposition}[theorem]{Proposition}
\theoremstyle{definition}

\theoremstyle{remark}
\newtheorem{remark}[theorem]{Remark}
\newtheorem{remarks}[theorem]{Remarks}
\newtheorem{example}[theorem]{Example}

\theoremstyle{fact}

\renewcommand{\Bbb}{\mathbb}

\def\norm#1{\left\Vert#1\right\Vert}

\def\N{{\Bbb N}}
\def\Z{{\Bbb Z}}
\def\S{{\mathbf S}}

\def\Homeo{{\mathrm{Homeo}}\,}

\def\QED{\nobreak\quad\ifmmode\roman{Q.E.D.}\else{\rm Q.E.D.}\fi}

\def\a{\alpha}

\def\implies{\Rightarrow}
\def\R{{\mathbb R}}

\def\o{\omega}

\def\sbs{\subset}

\def\e{\varepsilon}

\def\obr{^{-1}}
\def\stm{\setminus}

\def\Iso{\mathrm{Iso\,}}

\def\Int{\mathrm{Int\,}}

\newcommand{\card}{\rm{card\,}}

\def\mathbf{\Bbb}
\begin{document}

\title
[On metrizable enveloping semigroups] {On metrizable enveloping
semigroups}

%Authors
%    Information for first author
\author[E. Glasner]{Eli Glasner}
\address{Department of Mathematics,
Tel-Aviv University, Tel Aviv, Israel}
\email{glasner@math.tau.ac.il}
\urladdr{http://www.math.tau.ac.il/$^\sim$glasner}

%    Information for second author
\author[M. Megrelishvili]{Michael Megrelishvili}
\address{Department of Mathematics,
Bar-Ilan University, 52900 Ramat-Gan, Israel}
\email{megereli@math.biu.ac.il}
\urladdr{http://www.math.biu.ac.il/$^\sim$megereli}

%    Information for third author
\author[V.V. Uspenskij]{Vladimir V. Uspenskij}
\address{Department of Mathematics, 321 Morton Hall, Ohio
University, Athens, Ohio 45701, USA}
\email{uspensk@math.ohiou.edu}

\thanks{%
{\it 2000 Mathematics Subject Classification:}
%22F05 General theory of group and pseudogroup actions
%22A25 Representations of general topological groups and semigroups
%37B05 Transformations and group actions with special properties
%(minimality, distality, proximality, etc.)
%46B10 Duality and reflexivity
%46B22 Radon-Nikodym, Krein-Milman and related properties
%54H20 Topological dynamics
%54H15 Transformation groups and semigroups
Primary: 54H20.
Secondary 22A25, 22F05, 37B05, 46B10, 46B22, 54H15}

\date{September 20, 2006}

\keywords{almost equicontinuous, Asplund space,
%uu14s Asplund functions are no longer present!
Baire 1
function, barely continuous, enveloping semigroup, locally
equicontinuous, Radon-Nikod\'ym system, Rosenthal compact,
semigroup compactification,
sensitive dependence, weakly almost periodic,
 weak$^*$-topology.}

\begin{abstract}
When a topological group $G$ acts on a compact space $X$, its {\em
enveloping semigroup} $E(X)$ is the closure of the set of
$g$-translations, $g\in G$, in the compact space~$X^X$. Assume
that $X$ is metrizable. It has recently been shown by the first
two authors that the following conditions are equivalent: (1) $X$
%uu15s
%  I think that "hereditarily" (adverb) is better that "hereditary"
% (adjective) that we had in the next line; I corrected this
is hereditarily almost equicontinuous; (2) $X$ is hereditarily
non-sensitive; (3)
%mm10s
%for any $f\in C(X)$ the pseudometric $(x,y)\mapsto
%\sup\{|f(gx)-f(gy)|:g\in G\}$
for any compatible metric $d$ on $X$ the metric
$d_G(x,y):=\sup\{d(gx,gy): g\in G\}$
 defines a separable topology on~$X$; (4) the dynamical system
$(G,X)$ admits a proper representation on an Asplund Banach space.
We prove that these conditions are also equivalent to the
following: the enveloping semigroup $E(X)$ is metrizable.
\end{abstract}

\maketitle

\section{Introduction} \label{s:intro}

%mm10s
%By
A {\em dynamical system}, or a $G$-\emph{space},
%we mean
in this paper is a compact space $X$ (`compact' will mean `compact
and Hausdorff') on which a topological group $G$ acts
continuously. We denote such a system by $(G,X)$. For $g\in G$ the
{\em $g$-translation} (or {\em $g$-shift}) is the
self-homeomorphism $x\mapsto gx$ of $X$. If a nonempty subset $Y
\subset X$ is $G$-invariant, i.e. if $Y$ is closed under
$g$-shifts, then $Y$ is
%a $G$-\emph{subsystem} or
a $G$-\emph{subspace}.
The {\em enveloping semigroup} (or {\em Ellis semigroup}) of
$(G,X)$ is the closure of the set of $g$-shifts ($g\in G$) in the
compact space $X^X$, equipped with the product topology. Even for
simple dynamical systems on a compact metric space the enveloping
semigroup may be non-metrizable. For example, for the classical
Bernoulli shift (with $G:=\Z$) on the Cantor space
$X=\{0,1\}^{\Z}$, the enveloping semigroup $E(X)$ is homeomorphic
to $\beta \N$ (see \cite[Exercise 1.25]{Gl}). If $X$ is the unit
interval $[0,1]$ and $G=H_+[0,1]$ is the group of all
orientation-preserving homeomorphisms, then $E(X)$ is the
non-metrizable space of increasing and end-points-preserving
self-maps of $[0,1]$. If $X$ is a compact manifold without
boundary of dimension $>1$ and $G=\Homeo(X)$ is the group of all
self-homeomorphisms of $X$, then $E(X)$ is $X^X$.

On the other hand, if $G$ is an equicontinuous group of
homeomorphisms of a compact metric space $X$,
then $E(X)$ consists of continuous self-maps of $X$
and hence is metrizable. The same is true, more generally, if $(G,X)$ is
WAP (= Weakly Almost Periodic).
Recall that a function $f\in C(X)$ is {\em weakly almost periodic}
if its
$G$-orbit $\{^gf:g\in G\}$
lies in a weakly compact subset of the Banach space $C(X)$, and
$(G,X)$ is WAP if every $f\in C(X)$ is WAP. A dynamical system
$(G,X)$ is WAP if and only if $E(X)$ consists of continuous
self-maps of $X$ \cite{Ellis59, EllNer89}.

A generalization of WAP systems, called {\em Radon--Nikod\'ym} (RN
for short) systems, was studied in \cite{Me-NZ,GlasMegr}. To
define this notion, note that with every Banach space $V$ one can
associate a dynamical system $S_V=(H,Y)$ as follows: $H=\Iso(V)$
is the group of all linear isometries of $V$ onto itself, equipped
with pointwise convergence topology (or the compact-open topology,
the two topologies coincide on $H$), and $Y$ is the unit ball of
the dual space $V^*$, equipped with the weak$^*$-topology. The
action of $H$ on $Y$ is defined by $g\phi(v)=\phi(g\obr(v))$,
$g\in H$, $\phi\in Y$, $v\in V$. The continuity of this action can
be easily verified. A {\em representation} of a dynamical system
$(G,X)$ on a Banach space $V$ is a homomorphism of $(G,X)$ to
$S_V=(H,Y)$, that is, a pair of continuous maps $(h,\a)$, $h: G\to
\Iso(V)$ and $\a:X\to Y$, such that $h$ is a group homomorphism
and $\a(gx)=h(g)\a(x)$ for all $g\in G$ and $x\in X$. A
representation is {\em proper} if $\a$
%mm18s
%is an embedding.
is a topological embedding.
%

%mm18s breaking for an easier reading
A compact metric $G$-space $X$ is WAP if and only if $(G,X)$
admits a proper representation on a reflexive Banach space
%\cite[Corollary 6.10]{Me-NZ}, \cite[Theorem 7.6(1)]{GlasMegr}). A
%mm18s deleting ")"
\cite[Corollary 6.10]{Me-NZ}, \cite[Theorem 7.6(1)]{GlasMegr}. A
dynamical system is {\em Radon--Nikod\'ym} (RN) if it admits a
proper representation on an Asplund Banach space \cite[Definition
3.10]{Me-NZ}, \cite[Definition 7.5.2]{GlasMegr}. (If $G=\{1\}$, we
get the class of Radon--Nikod\'ym compact spaces in the sense of
Namioka \cite{N}.) Recall that a Banach space $V$ is {\em Asplund}
if for every separable subspace $E\sbs V$ the dual $E^*$ is
separable. Reflexive spaces and spaces of the
%uu15s
form
%type  %I think that "of the form" is better here than "of the type"
$c_0(\Gamma)$ are Asplund.
About the history and importance of Asplund spaces see for example
\cite{Bo,DGZ,F}.

Now assume that $X$ is a metrizable compact space. One
of the main results of \cite{GlasMegr} was a characterization of
RN-systems as those which are ``close to equicontinuous". To give
a precise statement we recall a few
definitions
from \cite{GW,AAB,Me-NZ,GlasMegr}.

Let $d$ be a compatible metric on $X$. We say that $(G,X)$ is {\em
non-sensitive} if for every $\e>0$ there exists a non-empty open
set $O\sbs X$ such that for every $g\in G$ the set $gO$ has
$d$-diameter $<\e$. (This property does not depend on the choice
of a compatible metric $d$.) A system $(G,X)$ is {\em hereditarily
non-sensitive} (HNS) if all
%mm10s
%subsystems (= $G$-invariant closed subsets)
closed $G$-subsystems
 are non-sensitive.

A system $(G,X)$ is {\em equicontinuous at $p\in X$} if for every
$\e>0$ there exists
a neighborhood $O$ of $p$ such that for every $x\in O$ and every
$g\in G$ we have
$d(gx,gp)<\e$. A system is {\em almost equicontinuous} (AE) if it is
equicontinuous at a dense
set of points, and {\em hereditarily almost equicontinuous} (HAE) if
every closed subsystem
is AE.

%uu14s
Denote by $Eq_\e$ the union of all open sets $O\sbs X$ such that
for every $g\in G$ the set $gO$ has diameter $<\e$. Then $Eq_\e$
is open and $G$-invariant. Let $Eq=\bigcap_{\e>0}Eq_\e$. Note that
a system $(G,X)$ is non-sensitive if and only if $Eq_\e\ne\emptyset$
for every $\e>0$, and $(G,X)$ is equicontinuous at $p\in X$ if and only
if $p\in Eq$. Suppose that $Eq_\e$ is dense for every $\e>0$.
Then $Eq$ is dense, in virtue of the Baire category theorem.
It follows that $(G,X)$ is AE.

If $(G,X)$ is non-sensitive and $x\in X$ is
%eee
%topologically transitive,
a {\em transitive point} ---
that is, $Gx$ is dense --- then for every $\e>0$ the open invariant set
$Eq_\e$ meets $Gx$ and hence contains~$Gx$. Thus $x\in Eq$.
If, in addition, $(G,X)$ is minimal (= all points are
%eee
%topologically
transitive), then $Eq=X$. Thus minimal non-sensitive systems are
equicontinuous (see \cite{AJ},
\cite[Theorem 1.3]{GW},
%eee
\cite{AAB96},
or \cite[Corollary 5.15]{GlasMegr}).

%mm10s
%A function $f\in C(X)$ is {\em Asplund} if the pseudometric $d_f$ on
%$X$ defined by
%$d_f(x,y)=\sup\{|f(gx)-f(gy)|:g\in G\}$ has the property that the
%pseudometric space
%$(X,d_f)$ is separable.

%mm18s we do not use it but perhaps it is worth to note that
%uu19s No, I think we better do not include this.
%The class of RN metric compact $G$-spaces is closed under closed
%$G$-subspaces, countable products and $G$-factors.

\begin{theorem}[{\cite[Theorem 9.14]{GlasMegr}}]
\label{th:GM}

%mm10s
%For a dynamical system $(G,X)$,
For a compact metric $G$-space $X$ the following conditions are
equivalent:
\begin{enumerate}
\item $X$ is RN;
\item $X$ is HNS;
\item $X$ is HAE;
\item every
%mm10s
nonempty
closed $G$-subspace $Y$ of $X$ has a point of
equicontinuity;
%mm10s
%\item every function $f\in C(X)$ is Asplund.
\item for any compatible metric $d$ on $X$ the metric
$d_G(x,y):=\sup_{g \in G} d(gx,gy)$ defines a separable topology
on $X$.
\end{enumerate}
\end{theorem}

It was proved in \cite{GlasMegr} that the equivalent conditions of
Theorem~\ref{th:GM}
imply that the enveloping semigroup $E(X)$ must be of cardinality
$\le 2^\omega$.
In fact, it was established in \cite[Theorem 14.8]{GlasMegr} that
$E(X)$ is Rosenthal compact
%uu19s
(see the first paragraph of Section~\ref{s:corol} for a definition),
and the question was posed whether this conclusion can be
strengthened to
``$E(X)$ is metrizable". This question was repeated in
\cite[Question 7.7]{Megr-OPIT2}.
The aim of the present paper is to answer this
question in the affirmative. Moreover, it turns out that
metrizablity of $E(X)$
%uu15s
% Eli, I still feel that something more emotional
% is appropriate here. I had "is just equivalent",
% you eliminated that, and  I agree that "just" is
% not the best word here. What about "in fact is equivalent"
% or "actually is equivalent"?
in fact
is equivalent to the conditions of
Theorem~\ref{th:GM}:

\begin{theorem}
\label{th:main}
Let $X$ be a compact metric $G$-space.
The following conditions are equivalent:
\begin{enumerate}
\item
the dynamical system $(G,X)$ is hereditarily almost equicontinuous
(HAE);
\item
the dynamical system $(G,X)$ is
RN, that is,
admits a proper representation on an Asplund
Banach space;
\item
the enveloping semigroup $E(X)$ is metrizable.
\end{enumerate}
\end{theorem}

After
%u7June recalling
providing
a few facts from general topology in
Section~\ref{s:rem}, we prove in
Section~\ref{s:=>}
the implication $(2)\implies(3)$ of Theorem~\ref{th:main},
in other words,
that for every RN compact metric $G$-space $X$ the
enveloping semigroup $E(X)$ is metrizable.
We prove the implication $(3)\implies(1)$ in Section~\ref{s:<=}.
Since $(1)$ and $(2)$ are known to be equivalent (Theorem~\ref{th:GM}),
this proves Theorem~\ref{th:main}. The implication $(1)\implies(3)$ is
thus proved via representations on Banach spaces; we give an alternative
direct proof in Section~\ref{s:alt}.
Some corollaries of the main theorem
are discussed in Section~\ref{s:corol}.

%e20s
We thank the anonymous referee
for a careful reading of the paper and for many
useful remarks.

\section{General topology: prerequisites}
\label{s:rem}

A subset of a topological space is {\em meagre} if it can be
covered by a countable family
of closed sets with empty interior. A space is
{\em Baire} if every meagre set has
empty interior, or, equivalently, if the intersection of any
countable family of dense
open sets is dense. Let us say that a (not necessarily continuous)
function $f:X\to Y$ is
{\em Baire~1} if the inverse image of every open set in $Y$ is
$F_\sigma$ (= the union
of countably many closed sets) in $X$.
%uu14s
According to this definition, Baire 1 functions need not be limits
of continuous functions. However, if the target space $Y$ is
metrizable (or, more generally, perfectly normal), then the limit
of every pointwise converging sequence of continuous functions is
Baire 1:

%mm12s one can organize it also as a Lemma
%uu14s Done!

\begin{proposition}[R. Baire]
\label{prop6}
If $Y$ is a metric space and $\{f_n : X \to Y\}$ is a sequence of continuous
functions converging pointwise to $f : X \to Y$ then f is Baire 1.
\end{proposition}

\begin{proof}
Let $U\sbs Y$ be open. There is a sequence $\{F_n\}$ of closed sets
such that $U=\bigcup F_n=\bigcup \Int F_n$, where Int denotes the interior.
Then $f\obr(U)$ is the union over $n$ and $k$ of the closed
sets $\bigcap_{i>n}f_i\obr(F_k)$.
\end{proof}

%It is easy to see that if $Y$ is a metric space and $f_n: X \to Y$
%is a sequence of continuous functions converging pointwise to $f:
%X \to Y$ then $f$ is Baire 1.
%

\begin{proposition}[R. Baire]
\label{prop1}
Let $f:X\to Y$ be Baire 1. If $X$ is Baire and $Y$
%mm10s
%has a countable base,
is separable and metrizable then
there exists a dense $G_\delta$-subset $A$ of $X$ such that $f$ is
continuous at every $x\in A$.
\end{proposition}

\begin{proof}
Let $\{U_n:n\in\o\}$ be a countable base for $Y$.
Write $f\obr(U_n)=\bigcup_k F_{nk}$, where
each $F_{nk}$ is closed, and consider the union $D$ of
the boundaries of all the $F_{nk}$'s.
Then $D$ is meagre, and it is easy to see that $f$ is continuous at
every point of the
dense $G_\delta$-set $A=X\stm D$.
\end{proof}

\begin{proposition}
\label{prop2} Let $f:X \to Y$ be a
%uu14s Misha, why did you delete the phrase
% "not necessarily continuous"? It makes no harm and helps
% avoid any confusion. I restored it. -V.U.
(not necessarily continuous)
function from a topological space $X$ to a separable metric space
$Y$. Suppose that the inverse image of every closed ball in $Y$ is
closed in $X$. Then $f$ is Baire 1.
\end{proposition}

\begin{proof}
Every open set $U$ in $Y$ is the union of a countable family of
closed balls, hence
$f\obr(U)$ is $F_\sigma$.
\end{proof}

We denote by $C(X,Y)$ the space of continuous maps from $X$ to
$Y$, equipped with the compact-open topology. If $X$ is compact
and $Y$ is metric, this topology is generated by the sup-metric.
%mm18s I am not sure that it is necessary
%but if you found it useful please remove comment
%uu19s OK, let us make this remark
If $X$ is compact metrizable then the group $\Homeo(X) \subset
C(X,X)$ of all self-homeomorphisms of $X$ is a separable and
metrizable topological group.

\begin{proposition}
\label{prop4}
%mm10s %uu14s
Let $X$ be Baire, $L$ separable metrizable, $K$ compact
metrizable, $Y$ dense in $K$.
If $f:X\to C(K,L)$ is a (not
necessarily continuous) function such that for every $y\in Y$
the function $x\mapsto f(x)(y)$
from $X$ to $L$ is continuous, then there exists
a dense $G_\delta$-subset $A$ of $X$ such that $f$ is continuous at
every $x\in A$.
\end{proposition}

%uu14s
The same result is true under the following assumptions: $Y=K$,
$K$ is compact but not necessarily metrizable, $X$ is regular
and strongly countably complete in the sense of Namioka~\cite{N2}.
For an easier proof of Namioka's theorem
that works under less restrictive assumptions,
see~\cite{SR}.

\begin{proof}
Equip $C=C(K,L)$ with the sup-metric
%mm10s
using a compatible metric $d$ on $L$.
Then $C$ is a separable metric space, and the inverse image under
$f$ of the closed ball of radius $r>0$ centered at $h \in C$ is
closed, being the intersection of the closed sets $\{x\in X:
d(f(x)(y), h(y))\le r\}$, $y\in Y$. Thus Propositions~\ref{prop1}
and~\ref{prop2} apply.
\end{proof}

A function $f:X\to Y$ is {\em barely continuous} if for every
closed non-empty $A\sbs X$ the restriction $f|A$ has a point of
continuity.
(This pun originates in a 1976 paper of E. Michael and I. Namioka,
\cite{MN}.)
%mm12s alternatively one can give also exact reference [Baire, 1899] if you have it
%uu14s
It is a classical fact
(contained in R.~Baire's Thesis, 1899)
that a function between Polish spaces is barely
continuous if and only if it is Baire 1 (see e.g. \cite[Theorem
24.15]{Kech}).
If $f:X\to Y$ is an onto barely continuous function between metric
spaces and $X$ is separable, then so is $Y$
%uu15s
\cite{MN}
(see also \cite[Lemma
6.5 and Proposition 6.7]{GlasMegr}).
%
%mm12s
%We'll
We will
need later a $G$-space version of this statement.

If $X$ and $Y$ are $G$-spaces, let us say that $f:X\to Y$ is
{\em $G$-barely continuous} if the restriction $f|A$ has a point of continuity
for every $G$-invariant closed non-empty subset $A\sbs X$. A {\em $G$-map}
between $G$-spaces is a map commuting with the action of~$G$.

\begin{proposition}
\label{prop5}
Let $X$ and $Y$ be metric spaces. Suppose that a (discrete) group $G$ acts
on $X$ by homeomorphisms and on $Y$ by isometries. Let $f:X\to Y$ be an onto
$G$-map. If $f$ is
$G$-barely continuous and $X$ is separable, then $Y$ is separable.
\end{proposition}
\begin{proof}
Pick $\e>0$.
Let $\a$ be the collection of all open subsets $U$ of $X$ such that $f(U)$
can be covered by countably many sets of diameter $\le\e$. Then $\a$ is
$G$-invariant and closed under countable unions. Since there exists a countable
subfamily $\beta\sbs \a$ such that $\cup \beta=\cup \a$, the family $\a$ has a largest
element, namely $V=\cup\a$. Let $A=X\stm V$. If $a\in A$ is a point of continuity
of
%mm18s
$f|A$,
%we write restrictions forms as $f|A$,  $A|C$...
%I would prefer $f|_C$,  $A|_C$  etc.
%your opinion ?
%uu19s I prefer $f|A$, but $f|_A$ also is acceptable.
there exists an open set $O\sbs X$ such that $a\in O$ and $f(O\cap
A)$ has diameter $\le \e$. Then $f(O\cup V)=f(O\cap A)\cup f(V)$
can be covered by countably many sets of diameter $\le \e$. Thus
$O\cup V\in \a$, in contradiction with the fact that $O$ meets the
complement of $\cup\a=V$. We have proved that $f|A$ has no points
of continuity. Since $A$ is closed and $G$-invariant, and $f$ is
$G$-barely continuous, it follows that $A$ is empty.

Thus $X=V\in \a$, and $Y$ can be covered by countably many sets of diameter $\le\e$.
Since $\e$ was arbitrary, $Y$ is separable.
\end{proof}

\begin{proposition}
\label{prop3}
The Banach dual $V^*$ of a non-separable Banach space $V$ is non-separable.
\end{proposition}

\begin{proof}
Construct a transfinite sequence $\{x_\a:\a<\o_1\}$ of unit
vectors in $V$ such that
for each $\a<\o_1$ %the distance from
the vector $x_\a$ does not belong to the closed linear
space $L_\a$ spanned
by the vectors $x_\beta$, $\beta<\a$. %, is $>1/2$.
For every $\a<\o_1$
find a functional
$f_\a\in V^*$ such that %$\norm{f_\a}=1$,
$f_a\in L_\a^\perp$
and $f_\a(x_\a)=1$. %$f_\a(x_\a)>1/2$.
All the pairwise distances between distinct $f_\a$'s are $\ge1$. %$>1/2$.
It follows that $V^*$,
considered with its norm topology, is not separable.
\end{proof}

%uu14s
%The following fact will be used in our proofs without an explicit
%reference:

\begin{proposition}
\label{prop7}
Let $f:X\to Y$ be a continuous onto map between compact spaces.
If $X$ is metrizable, then so is $Y$.
\end{proposition}

\begin{proof}
A compact space $K$ is metrizable if and only if it has a countable base
if and only if the Banach space $C(K)$ is separable. Note that $C(Y)$
is isometric to a subspace of $C(X)$ and hence is separable if $C(X)$ is
separable.

Alternatively, one can use Arhangelskii's theorem on coincidence
of the network weight and weight in compact spaces~\cite[Theorem
3.1.19]{Eng}.
%uu19s Misha, Eli has already checked the reference, in the 1989
%edition the theorem has the same number
This approach yields a stronger result: a compact space is
metrizable if it is the image under a continuous mapping of any
space with a countable base, compact or not.
\end{proof}

\section{Proof of Theorem~\ref{th:main}: Part 1}
\label{s:=>}

In this section we prove that for every RN compact metric
$G$-space $X$ the enveloping semigroup $E(X)$ is metrizable. Recall
that $X$ being RN means that $(G,X)$ has a proper representation
on an Asplund Banach space.

For a Banach space $V$ we denote by $S_V$ the dynamical system
%mm18s $Y$ was the phase space of $S_V=(H,Y)$ and not $X$ according to the Introduction
%so I replace $X$ by $Y$ in Proposition 3.1 (and $Y$ by $C$) in its proof
%$(\Iso(V), X)$, where $X$ is the unit ball of the dual space
$(\Iso(V), Y)$, where $Y$ is the unit ball of the dual space
$V^*$, equipped with the weak$^*$ topology.

We first prove the special case of the implication $(2)\implies(3)$
of Theorem~\ref{th:main},
when the dynamical system is of the form $S_V$, where $V$ is a Banach space
with a separable dual:

\begin{proposition}
\label{prop:=>}
Let $V$ be a Banach space with a separable dual,
$G=\Iso(V)$,
%mm18s
%$X$
$Y$
the compact
 unit ball of $V^*$ with the weak$^*$ topology, considered
as a $G$-space. Then the enveloping semigroup
%mm18s
%$E(X)$
$E(Y)$
is metrizable.
\end{proposition}
\begin{proof}
Let $K$ be the set of all linear operators of norm $\le 1$ on the
Banach space $V^*$. Consider the topology on $K$ inherited from
the product $(V^*)^{V^*}$, where each factor $V^*$ is equipped
with the weak$^*$ topology. Then $K$ is compact, being a closed
subset of the product $\prod_{f\in V^*} \norm{f}Y$. We claim that
$K$ is metrizable. Indeed, $V$ is separable
(Proposition~\ref{prop3}), hence $Y$ is metrizable, and so is each
ball $rY$, $r>0$. If $C$ is a norm-dense countable subset of $V^*$,
the restriction $A\to A|C$
defines a homeomorphism of $K$ onto a subspace of the product
%uu19s I changed $c\in C$ to $f\in C$. We used $f$ above,
%and we better avoid using $c$ and $C$, they are hard to distinguish
$\prod_{f\in C} \norm{f}Y$ of countably many metrizable compacta.
This proves our claim that $K$ is metrizable.

Restricting each operator $A\in K$ to $Y$, we obtain a
homeomorphism of $K$ with a compact subset $L$ of $Y^Y$. The
enveloping semigroup $E(Y)$ is the closure of the set $\{T^*|Y:
T\in G\}$ in $L$. Since $K$ is metrizable, so are $L$ and $E(Y)$.
\end{proof}

\begin{proposition} \label{weight}
Let $G$ be a separable topological group,
$X$ a compact metric $G$-space.
%uu19s
%Then
If $X$ is RN then
$(G,X)$ has a proper representation on a
% separable Asplund Banach space $V$.
Banach space with a separable dual.
\end{proposition}
\begin{proof}
%If we replace $G$ by its natural image
%in the topological group $\Homeo(X)$ (that is, the set of all
%$g$-shifts, $g\in G$), the enveloping semigroup will remain the
%same. Thus we may assume that $G$ is separable (and metrizable).
There exists a proper representation
$(h,\a):(G,X)\to S_V=(H,Y)$
for some Asplund $V$.
%We claim that $V$ can be chosen to have separable dual.
Since $\a(X)$ is metrizable, there exists a
countable subset $A\sbs V$ that separates points of $\a(X)$. Let
$W$ be the closed linear subspace of $V$ spanned by the union of
%uu19s
% $h(G)$-orbits is more precise but makes the notation heavier.
% There is no ambiguity about the action of $G$ on $V$. I changed
% $h(G)$-orbits to $G$-orbits, etc.
$G$-orbits of all points of $A$. Then $W$ is separable
%(because $h(G)$ is separable being a continuous image of $G$),
(note that the $G$-orbit of any point $v\in V$ is separable,
being a continuous image of $G$),
$G$-invariant, and the restriction map $V^*\to W^*$ is
one-to-one on $\a(X)$. It follows that $(G,X)$ admits a proper
representation on $W$. Since $V$ is Asplund and $W$ is separable,
the dual of $W$ is separable.
%Thus we may assume that $V=W$ and $V^*$ is separable.
\end{proof}

\begin{proposition} \label{part}
Let $X$ be a compact $G$-space. Suppose that $G_1$ is a
%uu19s
%topological
subgroup of $G$ and $X_1$ is a closed $G_1$-invariant subset of $X$.
%such that $X_1$ is $G_1$-invariant.
%uu19s
%with respect to the restricted action. Then
If $E(G,X)$ is metrizable then
%uu19s
% maybe just "then so is $E(G_1,X_1)$", rather than
% "$E(G_1,X_1)$ also is metrizable"?
$E(G_1,X_1)$
also is metrizable.
\end{proposition}
\begin{proof}
%uu19s
Consider the dynamical systems
$D=(G,X)$, $D_1=(G_1,X)$, and $D_2=(G_1,X_1)$.
The enveloping
semigroup $E(D_1)$ is a subspace of $E(D)$,
and there is a natural onto map
$E(D_1)\to E(D_2)$. If $E(D)$ is metrizable, then so are
$E(D_1)$ and $E(D_2)$ (Proposition~\ref{prop7}).
\end{proof}

We
now show that
for every RN compact metric $G$-space $X$
the enveloping semigroup $E(G,X)$ is
metrizable. 
%u20s
Since $E(G,X)$ depends only on the image of $G$ in $\Homeo(X)$,
we may assume that $G$ is a topological subgroup of $\Homeo(X)$
and hence separable.
By
Proposition \ref{weight} there exists a proper representation
$(h,\a):(G,X)\to S_V=(H,Y)$ for some
Banach space $V$ with a separable dual.
In virtue of Proposition~\ref{prop:=>}, the enveloping semigroup
of the system $S_V=(H,Y)$ is metrizable. Consider the dynamical
system $(h(G),\a(X))$. Its enveloping semigroup is metrizable by
Proposition \ref{part}. It remains to note that $E(h(G),\a(X))$
and $E(G,X)$ are isomorphic.

\section{Proof of Theorem~\ref{th:main}: Part 2}
\label{s:<=}
Let $X$ be a compact metric $G$-space such that $E(X)$ is
metrizable. We prove that
$X$ is HAE (= Hereditarily Almost Equicontinuous).

For every closed $G$-subsystem $Y$ of $X$ the enveloping semigroup
$E(Y)$ is metrizable, being
a continuous image
of $E(X)$. Thus it suffices to prove that $X$ is AE, that is,
that the system $(G,X)$
is equicontinuous at a dense set of points.

Consider the metric space $C=C(E,X)$ of all continuous maps from
$E=E(X)$ to~$X$, equipped with the sup-metric. For each $x\in X$
let $x^*\in C$ be the evaluation map defined by $x^*(e)=e(x)$,
$e\in E$. It is easy to see that the map $f:X\to C$ defined by
$f(x)=x^*$ is continuous at a point $x\in X$ if and only if
$(G,X)$ is equicontinuous at~$x$.
Thus we must prove that $f$ has a dense set
of points of continuity. This follows from
Proposition~\ref{prop4}, where $K=E$,
$L=X$ and $Y\sbs K$ is the set of all $G$-translations.

\section{An alternative proof of the implication $(1) \implies (3)$
in Theorem~\ref{th:main}}
\label{s:alt}
The implication $(1) \implies (3)$
in Theorem~\ref{th:main}:
%u7June
{\em if $X$ is metric and HAE, then $E(X)$ is metrizable}
-- was obtained in an indirect way, via representations on Banach spaces.
In this section we give a direct proof in the spirit of Section~\ref{s:<=}.

Consider the same evaluation map $f:X\to C(E,X)$ as in
Section~\ref{s:<=}. The assumption that $X$ is HAE
implies
that for every non-empty closed $G$-invariant subset $Y$ of $X$
the restriction $f|Y$ has a point of continuity. In other words,
$f$ is $G$-barely continuous in the sense of Section~\ref{s:rem}.

Consider the action of $G$ on $E$ given by $ge(x)=e(g\obr x)$
($g\in G$, $e\in E$, $x\in X$), and the action of $G$ on $C(E,X)$
given by $gh(e)=h(g\obr e)$ ($g\in G$, $h\in C(E,X)$, $e\in E$).
(We consider here $G$ as a group without topology; these actions need not
be continuous if $G$ is considered with its original topology.)
Then $G$ acts on $C(E,X)$ by isometries.
The evaluation map $f:X\to C(E,X)$
is a $G$-map.
Therefore we can apply Proposition~\ref{prop5}: $f(X)$ is a
separable subset of $C(E,X)$. Pick a dense countable subset $A$ of
$f(X)$. Since $f(X)$ separates points of $E$, so does~$A$. The
diagonal product $\bigtriangleup A:E\to X^A$ is therefore an
embedding. Since $X$ is metrizable and $A$ is countable, $X^A$ is
metrizable, and so is $E$.

\section{Some applications
%mm12s
and remarks}
 \label{s:corol}

\subsection{Tame dynamical systems}
\label{ss:tame}
For a topological space $X$ denote by $B_1(X)$ the space of all
Baire 1 real-valued functions on $X$, equipped with the pointwise
convergence topology.
A compact space $K$ is {\em Rosenthal} if it is homeomorphic to a subspace
of $B_1(X)$ for some Polish $X$.

In
\cite[Theorem 3.2]{GlasMegr} the following
dynamical
%\hfill\break
Bourgain-Fremlin-Talagrand dichotomy was established.

\begin{theorem}[A dynamical BFT dichotomy]\label{D-BFT}
Let $(G,X)$ be a metric dynamical system and let $E(X)$
be its enveloping semigroup. We have the following dichotomy.
Either
\begin{enumerate}
\item
$E(X)$ is separable Rosenthal compact, hence with cardinality
${\card}{E(X)} \leq 2^{\o}$; or
\item
the compact space $E$ contains a homeomorphic
copy of $\beta\N$, hence ${\card}{E(X)} = 2^{2^{\o}}$.
\end{enumerate}
\end{theorem}

In \cite{G} a dynamical system is called {\em tame}
if the first alternative occurs, i.e. $E(X)$ is
Rosenthal compact.
It is shown in \cite{G} that a minimal metrizable tame system
with a commutative acting group is PI.
(For the definition of PI and for more details on the structure theory
of minimal dynamical systems see e.g. \cite{G2000}.)
The authors of two recent works \cite{H} and \cite{KL}
improve
this result to show that under the same conditions the system
is in fact an almost 1-1 extension of an equicontinuous
system.

Under the stronger assumption that $E(X)$ is metrizable
Theorem \ref{th:main} now shows that the
commutativity assumption can be dropped and that the system is actually
equicontinuous. We get the following definitive result
in the spirit of R. Ellis' joint continuity theorem \cite{Ellis57}.

\begin{theorem} \label{t:equic}
A metric minimal system $(G,X)$ is equicontinuous if and only if its
enveloping semigroup $E(X)$ is metrizable.
\end{theorem}

\begin{proof}
It is well known
that the enveloping semigroup of a metric
equicontinuous system is a
metrizable compact topological group
(see e.g. \cite[Exercise 1.26]{Gl}).
Conversely, if $E(X)$ is metrizable then, by Theorem \ref{th:main},
$(G,X)$ is HAE and being also minimal it is equicontinuous
%uu14s
(see the paragraph before Theorem~\ref{th:GM}).
%
%(see \cite[Theorem 1.3]{GW} or \cite[Corollary 5.15]{GlasMegr}).
\end{proof}

%mm12s  goes to Remarks below
%\begin{remark}
%This answers negatively Problem 3.3 in \cite{G}.
%\end{remark}

%uAug25 (to the end of the subsection).

Our characterization of metrizable HNS systems as those having
metrizable enveloping semigroups should be compared with
the following theorem:

\begin{theorem}
A compact metric dynamical system $(G,X)$ is tame if and only if
every element of $E(X)$ is a Baire 1 function from $X$ to itself.
\end{theorem}

\begin{proof}
If $Y$ is a separable metric space and $B_1(X,Y)\sbs Y^X$ is the
space of Baire 1 functions from $X$ to $Y$, then every compact
subset of $B_1(X,Y)$ is Rosenthal. Indeed, $Y$ embeds in $\R^\N$,
hence $B_1(X,Y)$ embeds in $B_1(X,\R^\N)=B_1(X\times \N)$. In
particular, if $E(X)\sbs B_1(X,X)$, then $E(X)$ is Rosenthal,
which means that $(G,X)$ is tame.
Conversely, if $E(X)$ is
Rosenthal, then by the Bourgain-Fremlin-Talagrand theorem it is
Fr\'echet \cite{BFT}. (Recall that a topological space $K$ is {\em
Fr\'echet} if for every $A\subset K$ and every $x\in \overline{A}$
there exists a sequence of elements of $A$ which converges to
$x$.) In particular, every $p\in E(X)=\overline{G}$
%uu15s
(we may assume that $G\sbs \Homeo(X)$)
is the limit of a sequence
%of images
of elements of $G$
% in $C(X,X)$
and therefore of Baire class 1 (Proposition~\ref{prop6}).
\end{proof}

\begin{remarks}
%uu14s the "enumerate" environment does not look nice here,
% I enumerated the remarks by hand
%\begin{enumerate}
%    \item
    %We note that Theorem \ref{th:main} also resolves Problem 15.3 in \cite{GlasMegr}.
(1)   Note that Theorem \ref{th:main} resolves Problem 15.3 in \cite{GlasMegr}.
In fact, since the
%mm18s is it understood ? or we need to add [GW] reference (Colloq. Math.)
Glasner-Weiss examples
are metric and HNS
%mm18s
(see \cite[Section 11]{GlasMegr})
 we now know that their enveloping semigroups
are metrizable.
%mm12s
%\item

(2) Theorem \ref{t:equic} answers negatively Problem 3.3 in \cite{G}.

(3) In his paper \cite{Ellis93} Ellis, following Furstenberg's
classical work, investigates the projective action of $GL(n,\R)$
on the projective space $\mathbf{P}^{n-1}$. It follows from his
results that the corresponding enveloping semigroup is
not first countable.
In a later work \cite{Ak}, Akin studies the action of $G=GL(n,\R)$
on the sphere $\S^{n-1}$ and shows that here the enveloping
semigroup is first countable (but not metrizable).
The dynamical systems 
%u20s
$D_1=(G, \mathbf P^{n-1})$ and $D_2=(G, \S^{n-1})$ 
are tame but not RN. Note that $E(D_1)$ is Fr\'echet,
being a continuous image of a first countable space, namely $E(D_2)$.
\end{remarks}

\subsection{Distality and equicontinuity}

A dynamical system $(G,X)$ is {\em distal} if for any two distinct
points $x,y\in X$ the closure of the set $\{(gx,gy):g\in G\}$ in
$X^2$ is disjoint from the diagonal. If $X$ is metrizable and $d$ is
a compatible metric on $X$, this condition means that $\inf_{g\in
G}d(gx,gy)>0$.
Every equicontinuous system is distal.
By a theorem of Ellis a dynamical system $(G,X)$ is distal if and only if
its enveloping semigroup $E(X)$ is (algebraically) a group,
see \cite{Ellis57b}.
Note that this characterization implies that for any distal system
$(G,X)$ the phase space $X$ is the disjoint union of its minimal subsets.
In particular it follows that
{\em a point transitive distal system is minimal}.
(A dynamical system $(G,X)$ is {\em point transitive} if there
is some $x\in X$ for which the orbit $Gx$ is dense in $X$.)
As we have already observed, when $X$ is equicontinuous $E(X)$ is
actually a compact topological group.

One version of Ellis' famous joint continuity theorem says that a compact
dynamical system $(G,X)$ such that $E(X)$ is a group of continuous maps
is necessarily equicontinuous (see \cite{Ellis57}
and \cite[page 60]{Au}). Using Ellis's  characterizations of WAP and distality
this can be reformulated as follows:
{\em A distal WAP system is equicontinuous}.
We will now show that the WAP condition can not be much relaxed.

\begin{example}
The following is an example
of a dynamical system $(\mathbb{Z},X)$
which is distal, HAE, and its enveloping semigroup
$E(X)$ is a compact topological group isomorphic to the
2-adic integers. However, $(\mathbb{Z},X)$ is
not WAP and
{\em a fortiori} not equicontinuous.

Let ${\mathbf S = \R/\Z}$ (reals mod 1) be the circle.
Let $X = {\mathbf S}\times (\N \cup \{\infty\})$, where
$\N \cup \{\infty\}$ is the one point compactification of
the natural numbers.
Let  $T : X \to X$  be defined by:
$$
T(s,n) = (s + 2^{-n} , n), \quad      T(s,\infty)=(s,\infty).
$$
It is not hard to see that $E(X)$ is isomorphic to
%e18May
%$\Z_2$
%$\Z_2$ is universally understood as $\{0,1\}$.
the compact topological group $\mathbf{Z}_2$ of 2-adic integers.
The fact that $X$ is not WAP can be verified directly by observing that
$E(X)$ contains discontinuous maps.
Indeed, the map $f_a\in E(X)$ corresponding
to the 2-adic integer
$$
a=\dots10101=1+4+16+\dots
$$
can be described as follows: $f_a(s,n)=(s+a_n,n)$, where
$$
a_{2k}= \frac{2^{2k}-1}{3\cdot 2^{2k}}\to\frac{1}{3},\quad \quad
a_{2k+1}=\frac{2^{2k+2}-1}{3\cdot 2^{2k+1}}\to\frac{2}{3}.
$$
Geometrically this means that half of the circles are turned by
approximately $2\pi/3$, while the other half
are turned
by approximately the same angle in the opposite direction. The map $f_a$
is discontinuous at the points of the limit circle.
\end{example}

For a point transitive HAE system
distality is equivalent to equicontinuity
because,
as we have seen,
a distal point transitive system must be minimal
and a minimal HAE system is equicontinuous.

\subsection{Semigroup compactifications of groups}

A semigroup $S$ is {\em right topological} if it is equipped with
such a topology that for every $y\in S$ the map $x\mapsto xy$ from
$S$ to itself is continuous.
(Some authors use the term {\em left topological} for this.)
If for every $y\in S$ the self-maps
$x\mapsto xy$ and $x\mapsto yx$ of $S$ both are continuous, $S$ is
a {\em semitopological semigroup}.
A {\em right topological semigroup compactification} of a
topological group $G$ is a compact right topological semigroup $S$
together with a continuous semigroup morphism $G\to S$ with a
dense range
such that the induced action $G \times S \to S$ is continuous.
A typical example
is the enveloping semigroup $E(X)$ of a
dynamical system $(G,X)$
together with the natural map $G \to E(X)$.

{\em Semitopological semigroup compactifications} are defined
analogously.

We have the following direct corollaries of Theorem~\ref{th:main}.

\begin{corollary}\label{es-hae}
For a metric HAE system $(G,X)$ its
enveloping semigroup $E(X)$ is again a metrizable HAE system.
\end{corollary}

\begin{proof}
This follows from Theorem \ref{th:main} because the enveloping
semigroup of the flow $(G,E(G,X))$ is isomorphic to $E(G,X)$.
\end{proof}

\begin{corollary}\label{3-classes}
The following three classes coincide:
\begin{enumerate}
\item
Metrizable enveloping semigroups of %dynamical
$G$-systems.
\item
Enveloping semigroups of HAE metrizable %dynamical
$G$-systems.
\item
Metrizable right topological
semigroup compactifications of $G$.
\end{enumerate}
\end{corollary}
\begin{proof}
A dynamical system has the structure of a right topological
semigroup
compactification of $G$ if and only if it is the
enveloping semigroup of some dynamical system (see e.g.
\cite[Section 1.4]{Gl} or \cite[Section 2]{GlasMegr}).
\end{proof}

\begin{remark}
It is well known that the enveloping semigroup of a WAP dynamical
system is a semitopological semigroup compactification of $G$ (see
e.g. \cite[Section 1.4]{Gl} or \cite[Section 2]{GlasMegr}). Thus a
WAP version of Corollary~\ref{3-classes} (omitting part (1)) can
be obtained by changing `HAE' to `WAP' and `right topological
semigroup' to `semitopological semigroup'. Moreover, as was shown
in
\cite{Dow} (see also \cite[Theorem 1.48]{Gl}),
when the acting group $G$ is commutative, a point transitive WAP
system is always isomorphic to its enveloping semigroup, which in
this case is a commutative semitopological semigroup. Thus for
such $G$ the class of all metric, point transitive, WAP systems
coincides with that of all metrizable, commutative,
semitopological semigroup compactifications of~$G$.
\end{remark}

 %mm12s  restricting ourself by group case we certainly have a very nice framework
 %however we miss the important world of SEMIGROUP actions
\subsection{Semigroup actions}

Our main result (Theorem \ref{th:main}) remains true for semigroup
actions up to a more flexible version of HAE. Namely, we say that a
continuous action of a topological semigroup $S$ on a metric space
$(X,d)$ is HAE if for every (not necessarily $S$-invariant) closed
nonempty subset $Y$ there exists a dense subset $Y_0 \subset Y$
such that every point
$y_0 \in Y_0$
is a point of continuity of the
natural inclusion map $(Y,d|_Y) \to (X,d_S)$, where
$d_S(x,y):=\sup_{s\in S} d(sx,sy).$
(It is not hard to see that for $G$-group actions on compact metric
spaces this definition is equivalent to our old definition
which involved only $G$-invariant closed subsets.)
%
%is is equivalent (in our setting)
%to saying that $1_X: (X.d) \to (X,d_S)$ is fragmented
%$\pi_{\sharp}: Y \to C(S,X)$.
Then again HAE, RN and the metrizability of $E(X)$ are equivalent.
We omit the details.

%. That
%is it is enough to check the condition only for closed subsets $Y$
%that are $G$-invariant.
%

%indeed, if a $G$-space $X$ is HAE in the old sense then $(G,X)$ is HNS.
%equivalently, $1_X: X \to (X,d_G)$ is fragmented (see Lemma 9.4 in [GM]).
%equivalently, the inclusions $Y \to (X,d_G)$ are fragmented for every closed subset $Y \subset X$.
%Since $X$ (and hence, $Y$) is Polish and $(X,d_G)$ metric it follows that
%$Y \to (X,d_G)$ has a dense subset of points of continuity (see Prop. 6.6 in [GM]).
%this proof suggests also a semigroup version of HNS.

%mm12s  I give some details of the "generalized main thm".  Just for us.
%checking myself.
%\begin{theorem} \label{t:sem} Let $S$ be a topological semigroup and
%$X$ is a compact metric $S$-space. TFAE:
%\begin{enumerate}
%    \item $X$ is HAE.
%    \item $X$ RN.
%    \item $E(X)$ is is metrizable.
%\end{enumerate}
%\end{theorem}
%\begin{proof}
%
%(1) $\Rightarrow$ (2): Let $X$ be HAE. Then $1_X: (X.d) \to (X,d_S)$ is fragmented.
%Then $(X,d_S)$ is separable by \cite[Lemma 6.5]{GlasMegr}. Then
%for every $f \in C(X)$ the pseudometric space $(X,\rho_f)$ is also
%separable. This means that $f \in Asp(X)$. Thus, $C(X)=Asp(X)$.
%This means that $X$ is RN$_{app}$.
%But RN$_{app}$ = RN for compact metric spaces.
%
%(2) $\Rightarrow$ (3): Similar to the group case (as in the present paper).
%
%(3) $\Rightarrow$ (1): Similar to the group case (as in the present paper).
%\end{proof}

%I think that several results of [GM] can be extended to semigroup actions up to more flexible versions
%of HAE and HNS.
%


\begin{thebibliography}{10}

\bibitem{Ak}
E. Akin,
{\em Enveloping linear maps\/}, Topological dynamics and
applications, Contemporary Mathematics {\bfseries 215}, a volume
in honor of R.~Ellis, 1998, pp. 121-131.

\bibitem{AAB96}
E. Akin, J. Auslander, and K. Berg,
{\em When is a transitive map chaotic \/},
Convergence in Ergodic Theory and Probability,
Walter de Gruyter \& Co. 1996, pp. 25-40.

\bibitem{AAB}
E. Akin, J. Auslander, and K. Berg,
{\em Almost equicontinuity and
the enveloping semigroup\/}, Topological dynamics and
applications, Contemporary Mathematics  {\bfseries 215}, a volume
in honor of R.~Ellis, 1998, pp. 75-81.

\bibitem{Au}
J. Auslander,
{\em Minimal Flows and their Extensions\/},
Mathematics Studies
%uu19s Is 153 a volume number? We better make it bold: {\bf 153}
{\bf 153}, Notas de Matem\'atica, 1988.

\bibitem{AJ}
J. Auslander and J. Yorke,
{\em Interval maps, factors of maps, and chaos\/},
T\^ohoku Math. J. {\bf 32} (1980), 177-188.
%

%\bibitem{BJM}
% J.F. Berglund, H.D. Junghenn and P. Milnes,
%{\it Analysis on Semigroups}, Wiley, New York, 1989.

\bibitem{BFT}
J. Bourgain, D.H. Fremlin and M. Talagrand, {\it Pointwise compact
sets of Baire-measurable functions}, Amer. J. Math. {\bf 100}
(1978), 845-886.

\bibitem{Bo}
R.D. Bourgin,
\emph{Geometric Aspects of Convex Sets with the
Radon-Nikod\'ym Property,} Lecture Notes in Math., {\bf 993},
Springer-Verlag, 1983.

\bibitem{DGZ}
R. Deville, G. Godefroy and V. Zizler,
{\it Smoothness and renormings in Banach spaces}, Pitman
Monographs and Surveys in Pure and Applied Mathematics, {\bf 64},
Longman Scientific Technical, 1993.

\bibitem{Dow}
T. Downarowicz,
{\em Weakly almost periodic flows and hidden eigenvalues\/},
Topological dynamics and applications,
Contemporary Mathematics  {\bfseries 215},
a volume in honor of R. Ellis, 1998, pp. 101-120.

\bibitem{Ellis57}
R. Ellis,
{\em Locally compact transformation groups\/},
Duke Math.\  J.\ {\bfseries 24},
(1957), 119-126.

\bibitem{Ellis57b}
R. Ellis,
{\em Distal transformation groups\/},
Pacific J.\ Math.\ {\bfseries 8},
(1957), 401-405.

\bibitem{Ellis59}
R. Ellis,
{\it Equicontinuity and almost periodic functions,}
Proc. Amer. Math. Soc. {\bf 10} (1959), 637-643.

\bibitem{Ellis93}
R. Ellis,
{\it The enveloping semigroup of projective flows},
Ergod. Th. Dynam. Sys.
{\bf 13} (1993), 635-660.

\bibitem{EllNer89} R. Ellis and M. Nerurkar,
{\it Weakly almost periodic flows,}
Trans. Amer. Math. Soc. {\bf 313} (1989), 103-119.

\bibitem{Eng}
R. Engelking,
{\em  General topology\/},
revised and completed edition,
Heldermann Verlag, Berlin, 1989.

\bibitem{F}
M. Fabian,
{\em Gateaux differentiability of convex functions and
topology. Weak Asplund spaces\/},\ Canadian Math.\ Soc.\ Series of
Monographs and Advanced Texts, A Wiley-Interscience Publication,
New York, 1997.

\bibitem{G2000}
E. Glasner,
{\em Structure theory as a tool in topological dynamics\/},
Descriptive set theory and dynamical systems,
LMS Lecture note Series {\bfseries 277},
Cambridge University Press, Cambridge, 2000, 173-209.

\bibitem{Gl}
E. Glasner,
{\it Ergodic Theory via joinings}, Math. Surveys and
Monographs, AMS, {\bf 101}, 2003.

\bibitem{G}
E. Glasner,
{\it On tame dynamical systems}, Colloq. Math.
\textbf{105} (2006), 283-295.

\bibitem{GlasMegr} E. Glasner and M. Megrelishvili,
{\it Hereditarily non-sensitive
dynamical systems and linear representations},
Colloq. Math. \textbf{104} (2006), no. 2, 223-283.

\bibitem{GW}
E. Glasner and B. Weiss,
{\em Sensitive dependence on initial conditions\/},
Nonlinearity  {\bfseries 6}
(1993), 1067-1075.

%e20s
\bibitem{H}
W. Huang,
{\em Tame systems and scrambled pairs under an abelian group
action\/},
Ergod. Th. Dynam. Sys.
{\bf 26} (2006), 1549-1567.

\bibitem{Kech}
A.S. Kechris, {\em Classical descriptive set theory\/},
Springer-Verlag, Graduate texts in mathematics {\bfseries 156},
1991.

\bibitem{KL}
D. Kerr and H. Li,
{\em Independence in topological and $C^*$-dynamics\/},
to appear.

\bibitem{Me-NZ}
M. Megrelishvili,
{\em Fragmentability and representations of
flows\/}, Topology Proceedings {\bf 27} (2003), no.~2, 497-544.
See also: www.math.biu.ac.il/ ${\tilde{}}$ megereli.

\bibitem{Megr-OPIT2}
M. Megrelishvili,
{\it Topological transformation groups: selected
topics}, in book: Second edition of Open Problems in Topology
(Elliott Pearl, editor),
to appear. %Elsevier Science.

\bibitem{MN}
E. Michael and I. Namioka,
{\em Barely continuous functions},
%uu15s (English. Russian summary)
Bull. Acad. Polon. Sci. S\'er. Sci. Math. Astronom. Phys. {\bf 24} (1976),
no. 10, 889-892.

\bibitem{N2}
I. Namioka, {\em Separate continuity and joint continuity},
Pacif. J. Math. {\bf 51} (1974), 515-531.

\bibitem{N}
I. Namioka, {\em Radon-Nikod\'ym compact spaces and
fragmentability\/}, Mathematika {\bfseries 34} (1987), 258-281.

\bibitem{SR}
J. Saint Raymond, {\em Jeux topologiques et espaces de Namioka},
Proc. Amer. Math. Soc. {\bf 87} (1983), 499-504.
\end{thebibliography}
\end{document}